\documentclass[12pt]{amsart}

\def\CH{\mathcal {H}}
\def\bm{\mathbf{m}}
\def\bb{\mathbf{b}}
\def\ba{\mathbf{a}}
\def\bj{\mathbf{j}}

\def\bn{\mathbf{n}}
\def\CV{\mathcal {V}}
\def\CB{\mathcal {B}}

\def\CO{\mathcal {O}}

\def\CS{\mathcal {S}}
\def\H{\mathbb{H}}
\def\T{\mathbb{T}}
\def\Z{\mathbb{Z}}
\def\C{\mathbb{C}}
\def\N{\mathbb{N}}

\def\R{\mathbb{R}}

\def\1{\text{\bf {1}}}
\def\bs{\backslash}

\def\im{\mathop{\text{\rm{im}}}\nolimits}
\def\h{{\mathfrak h}}
\def\oline{\overline}
\def \la {\langle}
\def \ra {\rangle}

 \newtheorem{thm}{Theorem}[section]
 \newtheorem{cor}[thm]{Corollary}
 \newtheorem{lem}[thm]{Lemma}
 \newtheorem{prop}[thm]{Proposition}
 
 \newtheorem{rem}[thm]{Remark}

\newcommand{\be}{\begin{equation}}
\newcommand{\ee}{\end{equation}}
\newcommand{\bea}{\begin{eqnarray}}

\newcommand{\eea}{\end{eqnarray}}
\newcommand{\Bea}{\begin{eqnarray*}}
\newcommand{\Eea}{\end{eqnarray*}}

\newcounter{cnt1}
\newcounter{cnt2}
\newcounter{cnt3}
\newcommand{\blr}{\begin{list}{$($\roman{cnt1}$)$}
 {\usecounter{cnt1} \setlength{\topsep}{0pt}
 \setlength{\itemsep}{0pt}}}
\newcommand{\bla}{\begin{list}{$($\alph{cnt2}$)$}
 {\usecounter{cnt2} \setlength{\topsep}{0pt}
 \setlength{\itemsep}{0pt}}}
\newcommand{\bln}{\begin{list}{$($\arabic{cnt3}$)$}
 {\usecounter{cnt3} \setlength{\topsep}{0pt}
 \setlength{\itemsep}{0pt}}}
\newcommand{\el}{\end{list}}

\date{}
\begin{document}

\title[Heat kernel transform]
{Heat kernel transform for
nilmanifolds associated to the Heisenberg group
\vskip 1.5em {\tt By} }

\author[Kr\"otz, Thangavelu and Xu]{B.\ Kr\"otz, S.\ Thangavelu and Y.\ Xu}

\address{Max-Planck Institut f\"ur Mathematik\\
Vivatsgasse 7, \\
D-53111 Bonn {\it E-mail~:} {\tt
kroetz@mpim-bonn.mpg.de} }

\address{Department of Mathematics\\
Indian Institute of Science\\
Bangalore 560 012 {\it E-mail~:} {\tt
veluma@math.iisc.ernet.in} }

\address{Department of Mathematics\\
University of Oregon\\
Eugene, OR 97403 {\it E-mail~:} {\tt
yuan@math.uoregon.edu} }

\maketitle

\begin{abstract}
We study the heat kernel transform on a nilmanifold $ M $ of the Heisenberg
group. We show that the image of $ L^2(M) $ under this transform is a direct
sum of weighted Bergman spaces which are related to twisted Bergman and
Hermite-Bergman spaces.
\end{abstract}

\section{Introduction}
\setcounter{equation}{0}
Let us consider a complete analytic Riemannian 
manifold $M$ and let us denote 
by $k_t(x,y)$ the heat kernel on it. 
We fix $t>0$ and draw our attention to 
the map  
$$K: M\to L^2(M), \ \ m\mapsto k_t(m, \cdot)\, .$$  
This assignment is analytic and hence $K$ admits 
an analytic extension to a holomorphic map 

$$K^\sim: M_\C \to L^2(M), \ \ z\mapsto k_t^\sim (z,\cdot)\, $$
with $M_\C$ a Stein tube surrrounding $M$. 
Consequently,  we obtain a map 
$$T_t: L^2(M)\to \CO(M_\C); \ T_t(f)(z)=\int_M f(m) k_t^\sim(z,m)\ dm $$
which we call the the heat kernel transform. 
The basic problem now is to find appropriate tubes 
$M_\C$ and then characterize the image of $T_t$.  
This has been succesfully carried out for the following 
pairs $(M,M_\C)$:
\begin{itemize} 
\item $(\R^n, \C^n)$ \, \cite{B}. 
\item $(U, U_\C)$ with $U$ a compact Lie group 
and $U_\C$ its universal complexification \cite{H}. 
\item $(U/K, U_\C/ K_\C)$ with $U/K$ a compact 
symmetric space \cite{S}. 
\item $(\H, \H_\C)$ with $\H$ the Heiseberg group and 
$\H_\C$ its universal complexification \cite{KTX}.
\item $(G/K, \Xi)$ with $G/K$ a Riemannian symmetric space 
of the non-compact type and $\Xi$ the complex crown \cite{KOS}.
\end{itemize}
 
Let us mention that the image of $T_t$ can be very different in nature: a weighted 
Bergman space for $M=\R^n, U, U/K$,  a sum of 
two weighted Bergman spaces with oscillatory weight for $M=\H$, and, 
finally,  for  $M=G/K$ the image is not a Bergman space 
at all and needs to be characterized with tools from integral 
geometry. 
\par We observe,  that in all so far understood examples $M$ is
a symmetric space. Hence  one might ask 
if it is also possible 
to characterize $\im T_t$ for $M$ a locally 
symmetric space. This paper constitutes a modest 
step in that direction when we consider 
nilmanifolds associated to the Heisenberg group. 
\par To be precise, for $M=\Gamma\bs\H$ with 
$\Gamma<\H$ the standard lattice we show that 
$\im T_t$ is the sum of three weighted 
Bergman spaces two of which with oscillatory weight. 
To establish 
this theorem we lean on 
results for $\H$ \cite{KTX} as well as on 
deeper facts on the Hermite semigroup 
\cite{By}. Finally, let us mention 
that most of our methods extend to arbitrary 
discrete subgroups $\Gamma<\H$. 

\begin{center}
{\bf Acknowledgments}
\end{center}

The authors wish to thank the referee for making several useful comments on an
earlier version of the paper.

\section{Nilmanifolds associated to the Heisenberg group\\
and the heat kernel transform}
\setcounter{equation}{0}

\subsection{Nilmanifolds associated to the Heisenberg group}

Let us denote by $\H$ the $(2n+1)$-dimensional Heisenberg group. 
As a manifold $\H=\R^n \times \R^n \times \R$ and the group 
law is given by

$$(x,u,t) (x',u',t')=(x+x', u+u', t+t'+{1\over 2}
(u\cdot x'- x\cdot u'))\, .$$

Throughout $\Gamma$ will denote a discrete subgroup
of $\H$. With this data we form the nilmanifold $\Gamma\bs \H$. 
Sometimes we abbreviate $M =\Gamma\bs\H$. 
Here are some examples of $\Gamma$ we have in mind:

{\bf Examples:} (a) Let $\Gamma=\{0\}\times \{0\}\times\Z$. Then $\Gamma<\H$
is a discrete central subgroup. The quotient
$\H_{\rm red}=\Gamma\bs \H$ is the familiar reduced
Heisenberg group.
\par (b) The choice $\Gamma_{\rm st} =\Z^n\times\Z^n \times \frac12 \Z$ 
defines the standard lattice in $\H$.( This terminology is not standard!). 
The quotient
$\Gamma_{\rm st}\bs \H$ is a non-trivial circle bundle 
over the $2n$-torus  ${\T}^{2n}$ and hence compact. The fundamental group of
this compact manifold is the non-Abelian group $ \Gamma_{\rm st}.$ 
\par (c) The prescription $\Gamma=\Z^n\times \{0\}\times\{0\}$ 
defines an abelian discrete subgroup which does not intersect 
the center of $\H$.

For more about lattices in $\H$, especially their classification, 
we refer to \cite{Br}, \cite{F} and \cite{To}.
Write $N$ for the smallest connected subgroup containing $\Gamma$
($N$ coincides with the Zariski-closure of $\Gamma$ in 
the algebraic group $\H$ (cf. \cite{Ra}, Ch. II, Remark 2.6). )
We notice that $\Gamma$ becomes a lattice in the 
nilpotent group $N$.

\par Write ${\mathfrak n}$ for the Lie algebra of $N$. Then we find 
a subspace  ${\mathfrak v}\subset {\mathfrak \h}$ such that 
$\h={\mathfrak n}+{\mathfrak v}$. 
The prescription $V=\exp({\mathfrak v})$ defines a closed 
submanifold of $\H$ and the multiplication mapping 
$$N\times V\to \H, \ \ (n, v)\mapsto nv$$
is a polynomial homeomorphism. Consequently we obtain 
\begin{equation}\label{lattice}
\Gamma\bs \H\simeq \Gamma\bs N\times {\mathfrak v} 
\, ,\end{equation} 
where $\Gamma \bs N$ is compact and ${\mathfrak v}$ is a vector space. 

\par Write $dh$ for a Haar-measure on $\H$ which we 
normalize such that it coincides
with the Lebesgue measure once we identify $\H$ with $\R^{2n+1}$. 
Denote by $d(\Gamma h)$ the unique measure on $\Gamma\bs \H$ which satisfies
$$
\int_\H f(h)\,dh=\int_{\Gamma\bs \H} \sum_{\gamma\in \Gamma}  f(\gamma h) \,
d(\Gamma h)
$$
for all $f\in C_c(\H)$. The corresponding $L^p$-spaces shall be denoted by
$L^p(\Gamma\bs \H)$.

\par An important tool for us will be the averaging operator
$$A: C_c(\H)\to C_c(\Gamma\bs \H),  \ f\mapsto
A(f); \ A(f)(\Gamma h)=\sum_{\gamma\in \Gamma} f(\gamma h)$$
which is known to be continuous and onto. Observe that 
$A$ naturally extends to a surjective contraction 
$L^1(\H)\to L^1(\Gamma\bs \H)$, also denoted by $A$. Likewise 
(\ref{lattice}) implies that $A$ induces a continuous surjection of Schwartz 
spaces
$S(\H)\to S(\Gamma\bs \H)$. 

\subsection{Definition and basic properties of the heat kernel transform 
on $\Gamma\bs \H$}

The universal complexification of $\H_\C$ of 
$\H$ is simply $\H_\C=\C^n\times \C^n\times \C$ with 
holomorphically extended group law 
$$(z,w,\zeta)(z',w',\zeta')=(z+z', w+w', \zeta+ \zeta''+{1\over 2}
(w\cdot z'- z\cdot w'))\,.$$
Here $z\cdot w=\sum_{j=1}^n z_j w_j$.
If a real  analytic function $f$ on $\H$ admits
holomorphic extension to $\H_\C$, then $\tilde f$ shall 
be the notation for this so obtained function. 
  
\par We write $\Delta$ for the standard left 
Laplacian on $\H$ (cf. \cite{KTX}, Sect. 2.2) and denote 
by $k_t$ the corresponding heat kernel 
on $\H$. 
Explicitly 
$$k_t(x,u,\xi) = c_n \int_{-\infty}^\infty e^{-i\lambda \xi}e^{-t\lambda^2}
\left( \frac{\lambda}{\sinh(\lambda t)}\right)^n e^{-\frac{1}{4}\lambda
\coth(\lambda t)(x^2+u^2)} d\lambda\, $$
where $ c_n $ is a constant and $ x^2 = \sum_{j=1}^n x_j^2 $ , 
$ u^2 = \sum_{j=1}^n u_j^2 $ (cf. \cite{KTX}, Sect. 2.2). From the
above it is clear that $ k_t $ can be extended to a holomorphic function 
on $\H_\C$, namely 

\begin{equation}\label{eq=spec} k_t^\sim(z,w,\zeta) = c_n \int_{-\infty}^\infty e^{-i\lambda \zeta}
e^{-t\lambda^2}
\left( \frac{\lambda}{\sinh(\lambda t)}\right)^n e^{-\frac{1}{4}\lambda
\coth(\lambda t)(z^2+w^2)} d\lambda\, .\end{equation}
As before $ z^2 = \sum_{j=1}^n z_j^2 $ and $ w^2 = \sum_{j=1}^n w_j^2 $.
Note that $k_t^\sim$ is well behaved in the sense that 
its restriction to any bi-translate 
$g\H h\subset \H_\C$  with $h,g\in \H_\C$ is of rapid decay, 
i.e. $\H\ni x\mapsto k_t^\sim(gxh)$ is in $\CS(\H)$.   

\par For $t>0$ the {\it heat kernel transform} $T_t$ 
on $\H$ is defined by 

$$T_t : L^2 (\H)\to \CO(\H_\C), \ \ 
f\mapsto (f*k_t)^\sim$$
and one immediately verifies that 
\begin{itemize} 
\item $T_t$ is continuous (with $\CO(\H_\C)$ carrying the Fr\'echet
topology of compact convergence), 
\item $T_t$ is injective,  
\item $T_t$ is left-$\H$ equivariant 
\end{itemize}
(see \cite{KTX}, Sect. 3.1).

\par  It is not hard to see that $T_t$ extends to 
a map on tempered distributions 
$$\oline {T_t}: \CS'(\H)\to \CO(\H_\C), \ \ \nu\mapsto (\nu* k_t )^\sim $$
featuring the bulleted items from above. With that 
we turn to the heat kernel transform on the nil-manifold
$M=\Gamma\bs \H$. We often identify functions on $\Gamma\bs \H$ with 
$\Gamma$-invariant functions on $\H$. 
In this way we have $L^2(\Gamma\bs \H)\subset \CS'(\H)$ and 
with $T_t^\Gamma=\oline{T_t}|_{L^2(\Gamma\bs \H)}$ we obtain 
a continuous injection

$$T_t^\Gamma:  L^2(\Gamma\bs M) \to \CO(\Gamma\bs \H_\C)\, .$$

\begin{rem}  Recall that 
$\Delta$ was defined by the use of left-invariant
vector fields on $\H$ and so  factors to the Laplacian
$\Delta^\Gamma$ on the nilmanifold $\Gamma\bs \H$.
It is easy to see that 
$$ (e^{t\Delta^\Gamma}f)^\sim =T_t^\Gamma (f)$$
for all $f\in L^2(\Gamma \bs \H)$.
\end{rem}

It is useful to have an alternative way to describe 
$T_t^\Gamma$. 
Invoking the decomposition (\ref{lattice}) and the spectral 
resolution of $k_t^\sim$ (cf. (\ref{eq=spec}))  it is not hard to see 
that series 
$$K_t^\Gamma(\Gamma h, \Gamma z)=\sum_{\gamma\in\Gamma} k_t^\sim(h^{-1}\gamma w)
\qquad (h\in\H, z\in \H_\C)$$
converges uniformly in the first and locally uniformly in the second variable. 
As a result 
\begin{equation}
\label{growth} K_t^\Gamma\in \CS(\Gamma\bs \H_\C) \hat\otimes 
\CO(\Gamma\bs \H_\C)\, .  
\end{equation}
We observe 
\begin{equation}\label{heat1}
(T_t^\Gamma f)(\Gamma z) =\int_{\Gamma\bs \H} 
f(\Gamma h) K_t^\Gamma(\Gamma h, \Gamma z) \ d(\Gamma h)
\end{equation}
and deduce the inequality: 

\begin{lem} \label {Lemma B}
Let $Q\subseteq \H_\C$ be a compact subset. Then
$$C:= C(Q):=\sup_{z\in Q} \|K_t^\Gamma (\cdot, \Gamma z)\|_{L^2(\Gamma\bs \H)}
<\infty$$
and one has
$$
\sup_{z\in Q} |T_t^\Gamma(f)(\Gamma z)|\leq C \|f\|_{L^2(\Gamma\bs \H)}\, .$$
\end{lem}
We conclude this subsection with the averaging-equivariance. 

\begin{lem} \label{Lemma A}
Let $f\in \CS(\Gamma\bs \H)$ and $F\in \CS(\H)$
such that $A(F)=f$. Then
$$
T_t^\Gamma(f)=A(H_t(F))\, .
$$
\end{lem}

\begin{proof}
For all $z\in \H_\C$ we have
\begin{align*}
T_t^\Gamma(f)(\Gamma z) &= \int_\H f(h) k_t(h^{-1}z)\,dh \\
&=\sum_{\gamma\in \Gamma} \int_{\H} F(\gamma h) k_t(h^{-1}z) \,dh \\
&=\sum_{\gamma\in \Gamma} \int_{\H} F(h) k_t(h^{-1}\gamma z) \,dh \\
&=A(T_t(F))(\Gamma z)\  ,
\end{align*}
as was to be shown.
\end{proof}

\section{$\Gamma$-invariant distribution vectors and the Plancherel-Theorem for $L^2(M)$}
\setcounter{equation}{0}
Throughout this section and the next we confine ourselves to the standard lattice

$$\Gamma=\Gamma_{\rm st}=\Z^n \times \Z^n \times {1\over 2}\Z\, .$$

We will classify the $\Gamma$-invariant distribution vectors for 
the Schr\"odinger representation and relate this to the 
Plancherel decomposition $L^2(M)$. 
\par We wish to point out that the material collected below 
is all well known, see  \cite{GGPS}, \cite{A}, \cite{Br}, \cite{To}
and especially \cite{F} for a particularly nice treatment. The decomposition of
$ L^2(M) $ into irreducible pieces is due to Brezin \cite{Br} and we refer to
\cite{F} for an explicit exposition. In the terminology of \cite{F} the
lattice $ \Gamma_{\rm st} $ is isomorphic to $ \Gamma_l $ where $ l = (2,2,...,2) $ via the automorphism $ (x,u,t) \rightarrow (x,2u,2t) $.

In our exposition we adapt the  more general point of view of 
Gelfand et al. which, in our opinion, clarifies the underlying 
structure best.

\par To begin with we consider  a unimodular Lie group $G$ and let 
$\Gamma<G$ be a co-compact
lattice. Form $M=\Gamma\bs G$. One is interested in decomposing the right-regular 
representation $R$ on $L^2(M)$ into irreducibles. 
In this context one has a basic result (cf.\ \cite{GGPS}, Ch. 1, Sect. 4.6)

\begin{equation} \label{f1} L^2(M)\simeq \bigoplus_{\pi\in \hat G} m(\pi) \CH_\pi\, .
\end{equation}
Here, as usual,  $\hat G$ denotes the unitary dual of $G$ and the 
multiplicities $m(\pi)$ are all finite. By 
Frobenius reciprocity one has 
\begin{equation}\label{f2}  m(\pi)=\dim_\C (\CH_\pi^{-\infty})^\Gamma\, , \end{equation}
where $(\CH_\pi^{-\infty})^\Gamma$ denotes the space of $\Gamma$-invariant 
distribution vectors of $(\pi, \CH_\pi)$. 
Let us make the unitary equivalence in (\ref{f1}) explicit by writing down 
the intertwining operator. There exists
an inner product $\langle, \rangle_\pi$ 
on the finite dimensional $\C$-vectorspace ${\bf M}_\pi = 
(\CH_\pi^{-\infty})^\Gamma$ 
such that the map 

\begin{equation} \sum_{\pi\in\hat G}  
{\bf M}_\pi \otimes  \CH_\pi^\infty \to L^2(M), \ \ 
\sum \nu_\pi\otimes v_\pi\mapsto (\Gamma g \mapsto \sum_\pi 
\nu_\pi(\pi(g)v_\pi))
\end{equation}
extends to a unitary $G$-equivalence 
$$\bigoplus_{\pi\in\hat G}
{\bf M}_\pi\otimes  \CH_\pi\simeq L^2(M)\, .$$
\par Thus for our special situation $M=\Gamma\bs \H$ we have to determine 
two things: first ${\bf M}_\pi$,  and secondly
the inner product $\la\cdot
, \cdot \ra_\pi$ on this space. 
We turn to the details.

\par The first step in the decomposition of $ R $  
is the Fourier decomposition of $ f\in L^2(M) $ in the last variable. Note that $
F(x,u,\xi) $ is $\frac{1}{2}$- periodic in $ \xi $ and hence it has
the expansion
$$ F(x,u,\xi) = \sum_{k=-\infty}^\infty F^k(x,u)e^{4\pi ik \xi} $$ where
$ F^k(x,u) $ are the Fourier coefficients of $ F(x,u,\xi).$ Thus $ L^2(M) $
has the orthogonal direct sum decomposition
$$ L^2(M) = \sum_{k\in \Z}  \CH_k $$
where $ \CH_k $ is the set of all functions $ F \in L^2(M) $ satisfying
$ F(x,u,\xi) = e^{4\pi ik \xi} F(x,u,0).$
We now proceed to obtain further decomposition of  $ \CH_k $ for each
$ k \neq 0.$

Let $ \pi_\lambda, \lambda \in \R, \lambda \neq 0 $ be the Schr\"odinger
representations of $ \H $ realised on $ L^2(\R^n).$ Explicitly,
$$ \pi_\lambda(x,u,\xi)\varphi(v) = e^{i\lambda \xi}e^{i\lambda(x\cdot v
+\frac{1}{2}x\cdot u)}\varphi(v+u).$$ The subspaces $ \CH_k $ are invariant
under $ R $ and by Stone-von Neumann theorem $ R $ restricted to $ \CH_k $
decomposes into finitely many pieces each equivalent to $ \pi_{4\pi k}.$
For explicit decompositions of $ R $ we refer to \cite{Br},
\cite{To} and \cite{F}. There is no canonical way of
effecting the decomposition and here we get one such decomposition which
is suitable for our purpose of studying the heat kernel transform. For
some 'natural' decompositions of $ \CH_k $ we refer to Auslander and Brezin
\cite{A}.

\par As we already described earlier, the standard  way of constructing $ \Gamma $ invariant functions on $ \H $ is to 
start with
a tempered distribution $ \nu $ on $ \R^n $ which is $ \pi_\lambda(\Gamma) $ invariant
and consider $ F(x,u,\xi) = (\nu, \pi_\lambda(x,u,\xi)f) $ where $ f $ is a 
Schwartz function on $ \R^n.$ Let $ \nu $ be such a distribution; that is 
it verifies
$ (\nu, \pi_\lambda(h)f) = (\nu,f), ~ h \in \H.$ Then taking $ h = (0,0,j/2)
\in \Gamma, j \in \Z $ we are led to $ \pi_\lambda(h)f = e^{i\lambda j/2}f $ 
and $(\nu,e^{i\lambda j/2}f) = (\nu, f) $. This holds for all $ j \in \Z
$ if and only if $ \lambda = 4\pi k $ for some $ k \in \Z.$ Let us assume
$ k \neq 0 $ and write $ \rho_k = \pi_{4\pi k}.$

\begin{prop}Set  
$ A_k = \{ {\bj} \in \Z^n:
0 \leq j_1,j_2,....,j_n \leq 2k-1 \}.$ Then 
every tempered distribution $ \nu $ invariant under $ \rho_k(\Gamma) $ is
of the form
$ \nu = \sum_{{\bj} \in A_k} c_{\bj}\nu_{\bj} $ with $ \nu_{\bj}$ defined by
$$ (\nu_{\bj}, f) = \sum_{{\bm} \in \Z^n} \hat{f}(2k {\bm} +{\bj})\,\qquad (f\in \CS(\H)) $$ 
Here $ \hat{f} $ denotes the 
Fourier transform of the Schwartz class function $f$. 
\end{prop}
\begin{proof}  Since $ \lambda = 4k\pi $, $ e^{i\lambda \xi} = 1 $ for $ \xi \in
\frac{1}{2}\Z.$ The $ \rho_k(\Gamma) $-invariance of $ \nu $ shows that
$$ (\nu, \rho_k(0,{\bn},0)f) =(\nu,f(\cdot +{\bn})) = (\nu,f) $$ which means that
$ \nu $ is periodic. Consequently, $ \nu $ has the Fourier expansion
$$ (\nu,f) = \sum_{{\bn} \in \Z^n} c_{\bn} \hat{f}({\bn}) $$
where $ \hat{f}(\eta) = \int_{\R^n} f(x)e^{-2\pi i x\cdot \eta}\ dx.$ The
$ \rho_k(\Gamma) $ invariance of $ \nu $ also shows that
$$ (\nu, f) = (\nu,\rho_k({\bm},0,0)f) = (\nu,e^{4\pi k i {\bm} \cdot(\cdot)}f)$$
which translates into
$$ \sum_{{\bn} \in \Z^n} c_{\bn}\hat{f}({\bn}-2k{\bm}) =  \sum_{{\bn} \in \Z^n} c_{\bn}\hat{f}({\bn}).$$ This shows that $ c_{\bn} $ is a constant on the equivalence
classes in $ \Z^n/2k\Z^n.$ Thus
$$ (\nu,f) = \sum_{{\bj} \in A_k}c_{\bj} \sum_{\bm\in\Z^n} \hat{f}(2k{\bm} +{\bj}).$$ Defining
$ \nu_{\bj} $ accordingly we complete the proof.
\end{proof}

\begin{rem}\label{R1}  In view of the Poisson summation formula
$$ \sum_{{\bm }\in \Z^n} f(x+{\bm}) = \sum_{{\bm} \in \Z^n}\hat{f}({\bm})
e^{2\pi i
{\bm }\cdot x}\, , $$ valid for all functions $ f \in \CS(\R^n)$,  we have
$$ (\nu_{\bj},f) = (2k)^{-n}\sum_{{\bm} \in \Z^n}
e^{-\frac{\pi i}{k}{\bm} \cdot \bj} f(\frac{1}{2k}{\bm}).$$
\end{rem}
\begin{rem} Note that $\CS(\R^n)$ is the space 
of smooth vectors for the Schr\"odinger representation $\rho_k$, i.e.
$$L^2(\R^n)^\infty =\CS(\R^n)$$
in the standard representation theory terminology. 
Dualizing this identity we obtain 
$$L^2(\R^n)^{-\infty} =\CS'(\R^n)$$
and with it the Gelfand-triplet 
$$ L^2(\R^n)^\infty=\CS(\R^n)\hookrightarrow  L^2(\R^n)\hookrightarrow  \CS'(\R^n)=L^2(\R^n)^{-\infty}\, .$$ 
The above proposition then implies  that
$$\dim_\C (L^2(\R^n)^{-\infty})^\Gamma= (2k)^n$$
and, moreover, provides  an explicit basis for the space $(L^2(\R^n)^{-\infty})^\Gamma$. 
\par At this point it might be interesting to observe that there is 
an alternative way to construct elements of $(L^2(\R^n)^{-\infty})^\Gamma$, namely 
by averaging: Let $f\in\CS(\R^n)$. It is not difficult  to show that the 
series 
$$\sum_{\gamma\in\Gamma}\rho_k(\gamma)(f)$$
converges in $\CS'(\R^n)$ and defines a $\Gamma$-invariant 
element there. One establishes that the map 
$$S(\R)^n\to \CS'(\R^n)^\Gamma, \ \  f\mapsto \sum_{\gamma\in \Gamma} \rho_k(\gamma)(f)$$
is a continuous surjection.\end{rem}

At this point we determined the spectrum of $L^2(M)$, i.e. the occuring 
unitary irreducible representations, as well as the multiplicity 
space ${\bf M}_\pi$. In the sequel we abbreviate 
${\bf M}_k:={\bf M}_{\rho_k}$. As a last step we 
have to determine the unitary structure on ${\bf M}_k$ such that 
the map $ f \rightarrow ( \nu_{\bj},\rho_k(\cdot)f)$ becomes isometric. 
We already know that
$${\bf M}_k=\mathrm{span}\{ \nu_{\bj}\mid \bj\in A_k\}$$
and in the next step we want to show that the 
$\nu_{\bj}$ are in fact orthogonal. This is easy and follows from a
little group theory. In fact, let us define the finite 
group 

$${\bf F}_k:= (\Z/2k \Z)^n\, . $$ 
Then the prescription 
$$\Pi_k(x)(\nu):=\nu(\cdot + x)\qquad (\nu\in {\bf M}_k, x\in {\bf F}_k)$$
defines a representation of ${\bf F}_k$ on ${\bf M}_k$. 
Moreover it is clear that $\nu_{\bf j}$ is a basis 
of eigenvectors for this action; explicitly: 

\begin{equation} \label{f4}
\Pi_k(x)\nu_{\bj} =  e^{{\pi i\over k}  x\cdot {\bj}}\nu_{\bj} \qquad (\bj\in A_k)\, .\end{equation}
Futhermore for $f,g\in L^2(\R^n)$ and $\nu,\mu\in {\bf M}_k$ we set 
$F_\nu (x,u,\xi)=(\nu, \rho_k(x,u.\xi)f)$ and 
$G_\mu(x,u,\xi)=(\mu, \rho_k(x,u,\xi)g)$. Then one immediately 
verifies that 

$$(F_{\Pi_k(x)\nu}, G_{\mu})_{L^2(M)}= (F_\nu, G_{\Pi_k(-x)\mu})_{L^2(M)}\, .$$
Thus we conclude with (\ref{f4}) that $\nu_{\bj}$ is in fact 
an orthogonal basis (up to an uniform scalar) for ${\bf M}_k$. 
Furthermore it is indeed an orthonormal basis, as 
the next Lemma will show. 

\begin{lem}\label{Lem1} Let $f\in L^2(\R^n)$, 
${\bj}\in A_k$ and $F\in L^2(M)$ defined by the corresponding 
matrix coefficient 
$$F(x,u,\xi)=(\nu_{\bj}, \rho_k(x,u,\xi)f)\, .$$
Then 
$$\|F\|_{L^2(M)}=\sqrt{2}\cdot \|f\|_{L^2(\R^n)}\, .$$
\end{lem}

\begin{proof} The proof is a straightforward 
computation; we simply expand the terms: 

\begin{eqnarray*} F(x,u,\xi) & = &(\nu_{\bj}, \rho_k(x,u,\xi)f)\\
& = & (\nu_{\bj}, e^{4\pi ik \xi} e^{4\pi ik (x\cdot \cdot +{1\over 2}x\cdot u)}
f(\cdot + u)) \\
& = &e^{4\pi ik \xi} e^{2\pi ix\cdot u}  (\nu_{\bj}, e^{4\pi ik x\cdot \cdot}
f(\cdot + u)) \\
& = & {1\over (2k)^n} 
e^{4\pi ik \xi} e^{2\pi ix\cdot u} \sum_{\bm\in\Z^n} e^{-{\pi i\over k}
\bm\cdot \bj} e^{2\pi i x\cdot\bm}
 f\left({1\over 2k}\bm+u\right) \, .
\end{eqnarray*} 
In the last equation we used the characterization  of 
$\nu_j$ from Remark \ref{R1}. 
As $\Gamma\bs \H$ coincides with 
${\R/ {1\over 2}\Z}\times \R^n/\Z^n \times \R^n/\Z^n$
up to set of measure zero, 
we thus get 
\begin{eqnarray*}\|F\|_{L^2(M)}^2 & =& 
{2\over (2k)^{2n}} \int_{\R^n/\Z^n} \int_{\R^n/\Z^n}
\left|\sum_{\bm\in\Z^n} e^{-{\pi i\over k}
\bm\cdot \bj} e^{2\pi i x\cdot\bm}
 f\left({1\over 2k}\bm+u\right)\right|^2 \\ 
& & \qquad \ dx\ du \\ 
&=& {2\over (2k)^{2n}}
\int_{\R^n/\Z^n}
\sum_{\bm\in\Z^n}  \left|f\left({1\over 2k}\bm+u\right)\right|^2  \ du \\ 
 &=& {2\over (2k)^{n}}
\int_{\R^n/\Z^n}
\sum_{\bm\in\Z^n}  \left|f\left({1\over 2k}(\bm+u)\right)\right|^2  \ du\\ 
&=& {2\over (2k)^{n}}
\int_{\R^n}
\left|f\left({1\over 2k}u \right)\right|^2  \ du \\
&=& 2 \|f\|^2\ ,  
\end{eqnarray*}
as was to be shown. 
\end{proof}

Finally we make some remarks to the existing literature. 
We show that the matrix coefficients 
$F_{\bj}(x,u,\xi):=(\nu_{\bj},\rho_k(x,u,\xi)f)$ can be expressed as Weil-
Brezin (or Zak ) transforms studied in \cite{Br}. Consider the operator
$ V_k $ defined on the Schwartz class $ \CS(\R^n) $ by
$$ V_kf(x,u,\xi) = \sum_{{\bm} \in \Z^n} \rho_k(x,u,\xi)f({\bm}).$$ Written 
explicitly
$$ V_kf(x,u,\xi) = e^{4\pi ki\xi}e^{2\pi k i x\cdot u} \sum_{\bm \in \Z^n}
e^{4\pi k i {\bm} \cdot x} f(u+{\bm}).$$ It is easy to see that $ V_kf $ is
$ \Gamma $ invariant. For each $ \bj \in A_k $ we also define
$$ V_{k,\bj}f(x,u,\xi) = e^{2\pi i {\bj} \cdot x}V_kf(x,u,\xi).$$ These are
called the Weil-Brezin transforms in the literature.

\begin{prop}
\begin{enumerate}

\item  For each $ f \in \CS(\R^n) $ we have the relation
$$ F_{\bj}(x,u,\xi)= V_{k,\bj}g_{\bj}(u,-x,\xi) $$ where $ f $ 
and $ g_{\bj} $ are related by $ g_{\bj}(x) = \hat{f}(2kx+\bj).$
\item The transform $ V_{k,\bj} $ initially defined on $ \CS(\R^n) $ extends
to the whole of $ L^2(\R^n) $ as an isometry into $ \CH_k.$
\end{enumerate}
\end{prop}
 
\begin{proof} To prove the first assertion, 
a simple calculation shows that
$$ (\rho_k(x,u,\xi)f\hat{)}(s) = e^{4\pi ki \xi}e^{-2\pi ki x\cdot u}
e^{2\pi i s\cdot u}\hat{f}(s-2kx).$$ Hence it follows that $ F_{\bj}(x,u,\xi)
$ is given by
\bea
   e^{4\pi ki \xi}e^{-2\pi ki x\cdot u} 
\sum_{\bm \in \Z^n} e^{2\pi i(2k{\bm}+{\bj})\cdot u}\hat{f}(2k{\bm}+{\bj}-
2kx)
\eea
which simplifies to
$$ e^{4\pi ki \xi}e^{-2\pi ki x\cdot u}e^{2\pi i {\bj}\cdot u}
\sum_{\bm \in \Z^n} e^{4\pi i{\bm}\cdot u}\hat{f}(2k({\bm}-x)+{\bj}).$$
Setting $ g_{\bj}(s) = \hat{f}(2ks+\bj) $ and recalling the definition of
$ V_{k,\bj} $ we get $ F_{\bj}(x,u,\xi) = V_{k,\bj}g_{\bj}(u,-x,\xi).$
This shows that $ ( \nu_{\bj}, \rho_k(x,u,\xi)f) $ can also be defined on
the whole of $ L^2(\R^n).$ 
Finally, notice that the second assertion follows from the 
first one and our preceeding discussion. 
\end{proof}

Combining the two precceding propositions we get the following decomposition
of the spaces $ \CH_k.$ 

\begin{prop}
Let $ \CH_{k,\bj} $ be the span of functions of the form $F_{\bj}= 
(\nu_{\bj},\rho_k(\cdot)f) $
as $ f $ varies over $ L^2(\R^n).$ Then $ \CH_k $ is the orthogonal direct
sum of the spaces $\CH_{k,\bj}, {\bj} \in A_k.$
\end{prop}

\begin{rem} From the above proposition it follows that the restriction
of $ R $ to $ \CH_{k,\bj} $ is unitarily equivalent to $ \rho_k.$ The
intertwining operator is given by 
$$ U_{k,\bj}f(x,u,\xi) = (\nu_{\bj},\rho_k(x,u,\xi)f) $$ which is also equal
to the composition of $ V_{k,\bj} $ with the operators
$ f(s) \rightarrow \hat{f}(2ks+{\bj})$ and $ F(x,u,\xi) \rightarrow F(u,-x,\xi)
$. 
\end{rem}

An orthonormal basis for $ \CH_k $ can be obtained using the operators $
U_{k,\bj}.$ Let $ \Phi_\alpha, \alpha \in \N^n $ be the normalised Hermite
functions on $ \R^n.$ Then the functions $ U_{k,\bj}\Phi_\alpha(x,u,\xi) $ form
an orthogonal system in $ \CH_{k,\bj}.$ With suitable normalising constants
$ c_{\alpha,\bj} $ the functions $c_{\alpha,\bj}U_{k,\bj}\Phi_\alpha $ as
$ \alpha $ ranges over $ \N^n $ and $ {\bj} \in A_k $ form an orthonrmal basis
for $ \CH_k.$

\section{The image of the heat kernel transform}
\setcounter{equation}{0}
In this section we determine the image of $T_t^\Gamma$ for 
$\Gamma=\Gamma_{\rm st}$ the standard lattice. 
To simplify notation we often write $T_t$ instead of $T_t^\Gamma$ and drop 
the $\sim$ for the holomorphic extension of a function.

\par As $ L^2(M) $ is the direct sum of $ \CH_k $ as $ k $
ranges over all integers, the image of $ L^2(M) $ under $ T_t $ will be
the direct sum of $ T_t(\CH_k), $ the image of $ \CH_k $ under $ T_t.$
We first settle the case $ k =0.$ Recall that functions $ f \in \CH_0 $
are independent of $ \xi $ and hence we think of them as functions on the
$2n$-torus $ \T^n \times \T^n.$  An easy calculation (use (\ref{eq=spec}))
shows that the
function $ f*k_t $ is given by the ordinary convolution
$$ f*k_t(z,w) =c_n t^{-n} \int_{\R^n}\int_{\R^n} f(x',u')
e^{-\frac{1}{4t}((z-x')^2+(w-u')^2)} \ dx'\ du'.$$ Note that $ f*k_t $ is
an entire function on $ \C^n \times \C^n $ which satisfies $ f*k_t(z+
{\bm},w+{\bn}) = f*k_t(z,w) $ for all $ {\bm},{\bn} \in \Z^n.$ Thus the
heat kernel transform when restricted to $ \CH_0 $ is nothing but the heat
kernel transform on the torus $ \T^n \times \T^n $ and the image has been
characterised.

\begin{thm} An entire function $ F(z,w) $
of $ 2n$ complex variables  belongs to $T_t(\CH_0) $ if and only if
 $ F(z+{\bm},w+{\bn}) = F(z,w) $ for 
all $ {\bm},{\bn} \in \Z^n $ and
$$\|F\|^2=\int_{\R^{2n}} \int_{[0,1)^{2n}}|F(z,w)|^2 e^{-\frac{1}{2t}(y^2+v^2)}
\ dx\ du \ dy\ dv < \infty.$$
Moreover $\|T_t(f)\|=\|f\|$ for all $f\in \CH_0$. 
\end{thm}

Thus the members of $T_t(\CH_0) $ are precisely the functions from the
classical weighted Bergman space associated to the standard Laplacian
on $ \R^{2n}$ that are periodic in the real parts of the variables. 

We now consider the image of $ \CH_k$ for $k\neq 0$. For the description of $ T_t(\CH_k) $ 
we need to recall
several facts about Hermite-Bergman spaces and twisted Bergman spaces. Given
a nonzero $\lambda \in \R $ consider the kernel
$$ p_t^\lambda(x,u) = 
c_n \lambda^n(\sinh(\lambda t))^{-n}e^{-\frac{1}{4}\lambda \coth(\lambda t)(x^2+u^2)}.
$$
This kernel is related to $ k_t $ via
$$ k_t^\lambda(x,u) = e^{-t\lambda^2} p_t^\lambda(x,u) $$ where for a
function $ f(x,u,\xi) $ on the Heisenberg group we use the notation
$$ f^\lambda(x,u) = \int_{-\infty}^\infty f(x,u,\xi)e^{i\lambda \xi} \ d\xi.$$
Given a function $ f \in L^2(\R^{2n}) $ the $\lambda$-twisted convolution
$$ f*_\lambda p_t^\lambda(x,u) = \int_{\R^{2n}}f(x',u')p_t^\lambda(x-x',u-u')
e^{-i\frac{\lambda}{2}(u\cdot x'-x\cdot u')}\ dx'\ du' $$ 
extends to $ \C^{2n} $ as an entire function. This transform is called the
twisted heat kernel transform and in [7] we have studied the image
of $ L^2(\R^{2n}) $ under this transform.

The image turns out to be the twisted Bergman space $ \CB_t^\lambda $ which
is defined as follows. An entire function $ F(z,w) $ on $ \C^{2n} $ belongs to
$\CB_t^\lambda $ if and only if
$$ \int_{\C^n}\int_{\C^n}|F(z,w)|^2 W_t^\lambda(z,w) \ dz \ dw <\infty $$
where
$$ W_t^\lambda(x+iy,u+iv) = e^{\lambda(u\cdot y-v\cdot x)}p_{2t}^\lambda(
2y,2v).$$
In [7] it  has been shown that the image of $ L^2(\R^{2n}) $ under the twisted
heat kernel transform is precisely $ \CB_t^\lambda.$ 

The connection between the twisted heat kernel trasnform and the heat kernel 
transform on the nilmanifold $ M $ is the following. If $ F \in \CH_k $ then
an easy calculation shows that
\bea
 F*k_t(x,u,\xi) = e^{-t(4\pi k)^2} e^{4\pi ik \xi} G*_{-4\pi k}p_t^{-4\pi k}
(x,u)
\eea
where $ G(x,u) = F(x,u,0).$ Thus we are led to consider $\lambda$-twisted 
convolution with $ p_t^\lambda.$ Observe that  when $ F \in \CH_k $ the
function $ G(x,u) = F(x,u,0) $ satisfies
$$ G(x+{\bm},u+{\bn}) = e^{2\pi ik(u\cdot {\bm}-x\cdot {\bn})}G(x,u).$$ 
Thus the
entire extension of $  G*_{-4\pi k}p_t^{-4\pi k} $ also inherits the same
property and  we expect the image of $ \CH_k $ under the heat kernel
transform to consist of entire functions of the form
\begin{equation} \label{EE} F(z,w,\zeta) = e^{4\pi ik \zeta} G(z,w)\, ,\end{equation}
where $ G $ has the above 
transformation property under translation by $ \Z^n \times \Z^n.$ We define
$ \CB_{t,\Gamma}^{4\pi k} $ to be the space of all entire functions $ G(z,w)
$ having the transformation property
$$G(z+{\bm},w+{\bn}) = e^{2\pi ik(w\cdot {\bm}-z\cdot {\bn})}G(z,w)$$
which are square integrable with respect to $ W_t^{-4\pi k} $, i.e.
$$ \int_{\R^{2n}}\int_{[0,1)^{2n}}|G(z,w)|^2 W_t^{-4\pi k}(z,w)\ dx\ du
\ dy\ dv < \infty\, .$$ 

\subsection{Diagonalization of $ \CB_{t,\Gamma}^{4\pi k} $}
In this subsection we show
that $ \CB_{t,\Gamma}^{4\pi k} $ admits a natural symmetry 
of the finite group ${\bf F}_k$. 
To begin with we note that the prescription 
\begin{equation} \label{A1} (\tilde\Pi_k(x) G)(z,w)=e^{-i\pi x\cdot w} G\left(z+{x\over 2k}, w\right)\qquad 
(x\in {\bf F}_k, G\in \CB_{t,\Gamma}^{4\pi k})\, \end{equation}
defines an action of ${\bf F}_k$ on  $ \CB_{t,\Gamma}^{4\pi k}$. 
Moreover, 

\begin{lem} The representation $(\tilde \Pi_k,  \CB_{t,\Gamma}^{4\pi k})$ of 
${\bf F}_k$ is unitary. 
\end{lem}

\begin{proof}  We have to show that 
$$\|\tilde \Pi_k(s)G\|^2 =\|G\|^2$$
for all $G\in \CB_{t,\Gamma}^{4\pi k}$ and  $s\in {\bf F}_k$. 
The verification is a straightforward computation:

\begin{align*} \|\tilde \Pi_k(s)G\|^2 &=
\int_{\R^{2n}}\int_{[0,1)^{2n}}\left|(\tilde\Pi_k(s)G)(x+iy,u+iv)\right|^2 \cdot \\
&\qquad W_t^{-4\pi k}(x+iy,u+iv)\ dx\ du\ dy\ dv\\
&=\int_{\R^{2n}}\int_{[0,1)^{2n}}\left|(\tilde\Pi_k(s)G)(x+iy,u+iv)\right|^2 \cdot\\
&\qquad e^{-4\pi k (u\cdot y- v\cdot x)}
p_{2t}^{-4\pi k}(2y,2v)\ dx\ du\ dy\ dv \\
&=\int_{\R^{2n}}\int_{[0,1)^{2n}}\left|G(x+s/2k+iy,u+iv)\right|^2 \cdot\\
&\qquad e^{2\pi s\cdot v}   e^{-4\pi k (u\cdot y- v\cdot x)}
p_{2t}^{-4\pi k}(2y,2v)\ dx\ du\ dy\ dv\\
&=\int_{\R^{2n}}\int_{[0,1)^{2n}}\left|G(x+iy,u+iv)\right|^2 \cdot\\
&\qquad e^{2\pi s\cdot v}   e^{-4\pi k (u\cdot y- v\cdot (x-s/2k))}
p_{2t}^{-4\pi k}(2y,2v)\ dx\ du\ dy\ dv\\
&=\int_{\R^{2n}}\int_{[0,1)^{2n}}\left|G(x+iy,u+iv)\right|^2 \cdot\\
&\qquad e^{-4\pi k (u\cdot y- v\cdot x)}
p_{2t}^{-4\pi k}(2y,2v)\ dx\ du\ dy\ dv\\
&=\|G\|^2\, .
\end{align*}
\end{proof} 

In the sequel we often abbreviate and write $\lambda$ for $4\pi k$. 
Define 
for each $ {\bj} \in A_k $ a subspace of $ \CB_{t,\bj}^\lambda $ 
of  $ \CB_{t,\Gamma}^\lambda $ as follows: $G\in  \CB_{t,\Gamma}^\lambda$ 
belongs to  $ \CB_{t,\bj}^\lambda $ if and only if 
$$ G(z+\frac{1}{2k}{\bm},w) = e^{\pi i {\bm}\cdot(w+\frac{1}{k}{\bj})}
G(z,w)\, . $$
The previous lemma then implies:

\begin{cor}\label{corr} $\CB_{t,\Gamma}^\lambda $ is the orthogonal direct sum of 
$\CB_{t,\bj}^\lambda, {\bj} \in A_k.$
\end{cor}

\subsection{Characterization of $ \CB_{t,\Gamma}^\lambda $}
The aim of this subsection 
is to prove that $T_t$ maps $\CH_k$ onto 
$\CB_{t,\Gamma}^\lambda$ (with $\CB_{t,\Gamma}^\lambda$ interpreted
as a subspace of $\CO(\Gamma\bs \H_\C)$ via (\ref{EE})). 
Let us begin with the easy half, the isometry of the map. 
For that we have to introduce a useful 
technical tool, namely twisted averages. 

Recall  that $\CH_k^\infty$ is the 
space of all functions $F$ in
$ \CS(\Gamma\bs\H)$ which satisy 
\begin{equation}\label{GG}  F(x,u,\xi) = e^{-i\lambda \xi}F(x,u,0)\, .\end{equation} 
Often it is convenient to identify 
functions $f$ on $\R^{2n}$ with functions 
$F$ on $\H$ which transform as (\ref{GG}) via

\begin{equation} \label{HH}
f\leftrightarrow F , \ \  F(x,u,\xi)= e^{-i\lambda \xi} f(x,u)\, .\end{equation}
With this terminology we record 
for $F\in \CH_k$

\begin{equation}\label{FF} T_t(F)(z,w,\zeta)=  e^{i\lambda \zeta} (f*_\lambda p_t^\lambda)(z,w)\, . 
\end{equation}

For $f\in \CS(\R^{2n})$ we define its twisted average 
by 
$$ A_\lambda f(x,u) = e^{i\lambda \xi} \sum_{\gamma \in \Gamma} F(\gamma
(x,u,\xi)).$$ More explicitly,
\bea
 A_\lambda f(x,u) = \sum_{(\ba,\bb,0)\in \Gamma} e^{i\frac{\lambda}{2}
(u\cdot \ba- x\cdot \bb)}f(x+\ba,u+\bb).
\eea
We note that $A_\lambda$ maps $\CS(\R^{2n})$ 
onto $\CH_k^\infty$ (modulo the identification (\ref{HH})). 
Further we note that 

\begin{equation}\label{KK}
(A_\lambda f) *_\lambda p_t^\lambda= A_\lambda(f*_\lambda p_t^\lambda)
\end{equation}

\begin{prop}\label{stand} For all $f\in \CS(\R^{2n})$ one has
$$\|T_t(A_\lambda (f))\|_{\CB_{t,\Gamma}^\lambda}^2=\|A_\lambda(f)\|_{L^2(M)}^2\, .$$
 In particular, the map 
$$T_t: \CH_k\to \CB_{t,\Gamma}^\lambda, \ \ F \mapsto f*_\lambda p_t^\lambda$$
is isometric. 
\end{prop}

\begin{proof}
Let $f\in \CS(\R^{2n})$. Then, by  (\ref{KK}), 
$$ (A_\lambda f)*_\lambda p_t^\lambda(z,w) = \sum_{(\ba,\bb,0)\in \Gamma} 
e^{i\frac{\lambda}{2}
(w\cdot \ba- z\cdot \bb)}f*_\lambda p_t^\lambda(z+\ba,w+\bb).$$ 
We can obtain pointwise estimates for the function 
$ f*_\lambda p_t^\lambda(z+\ba,w+\bb)$ which shows that the above series
actually converges. Therefore,
\begin{align*}
\|T_t(A_\lambda(f))\|^2 & = \int_{\R^{2n}}
 \int_{\R^{2n}/\Z^{2n}} |(A_\lambda f)*_\lambda p_t^\lambda(z,w)|^2 W_t^\lambda(z,w) \ dz \ dw \\
& = \int_{\C^{2n}} 
f*_\lambda p_t^\lambda(z,w)
\overline{(A_\lambda f)*_\lambda p_t^\lambda(z,w)}W_t^\lambda(z,w)\ dz \ dw 
\end{align*}
and we used the transformation property of the weight function 
$W_t^\lambda$. 
Further expansion yields 

\begin{align*}
\|T_t(A_\lambda(f))\|^2 & 
 = \sum_{(\ba,\bb,0)\in \Gamma}\int_{\C^{2n}}
f*_\lambda p_t^\lambda(z,w)\cdot \\
&\quad 
e^{-i\frac{\lambda}{2}
(w\cdot \ba- z\cdot \bb)}\overline{f*_\lambda p_t^\lambda}(z+\ba,w+\bb)
W_t^\lambda(z,w)\ dz \ dw\ .
\end{align*}
We recall that 
$ W_t^\lambda $ is
the weight function for the twisted Bergman space $ \CB_t^\lambda $ (see \cite{KTX})
and 
obtain further

\begin{align*}
\|T_t(A_\lambda(f))\|^2 & 
 = \sum_{(\ba,\bb,0)\in \Gamma}\int_{\R^{2n}}
f(x,u)
e^{-i\frac{\lambda}{2}(u\cdot \ba- x\cdot \bb)}\overline{f(x+\ba,u+\bb)}\ dx\ du\\
&  = \int_{\R^{2n}/ \Z^{2n}} |A_\lambda f(x,u)|^2 \ dx \ du\\
&= \|A_\lambda f\|^2\, ,
\end{align*}
which completes the proof. 
\end{proof}

We turn to the more difficult part, namely that 
$T_t$ maps $\CH_k$ onto $\CB_{t,\Gamma}^\lambda$. 
This 
will be proved by establishing a connection between twisted
Bergman spaces and Hermite-Bergman spaces which we proceed to describe now.
For each nonzero $ \lambda \in \R $ let us consider the scaled Hermite operator
$ H(\lambda) = -\Delta +\lambda^2 |x|^2 $ on $ \R^n $ whose eigenfunctions
are provided by the Hermite functions
$$ \Phi_\alpha^\lambda(x) =|\lambda|^{\frac{n}{4}}\Phi_\alpha(\sqrt{|\lambda|}
x), x \in \R^n, \alpha \in \N^n.$$ The operator $ H(\lambda) $ generates the
Hermite semigroup $ e^{-tH(\lambda)} $ whose kernel is explicitly given by
$$ K_t^\lambda(x,u) = \sum_{\alpha \in \N^n}e^{-(2|\alpha|+n)|\lambda|t}
\Phi_\alpha^\lambda(x)\Phi_\alpha^\lambda(u).$$ Using Mehler's formula
(see \cite{Th}) the above series can be summed to get
\bea\label{Mehler}
 K_t^\lambda(x,u) = c_n (\sinh(\lambda t))^{-\frac{n}{2}}
(\cosh(\lambda t))^{-\frac{n}{2}}
\eea
$$
\times e^{-\frac{\lambda}{4}\tanh(\lambda t)(x+u)^2}
e^{-\frac{\lambda}{4}\coth(\lambda t)(x-u)^2}.
$$
The image of $ L^2(\R^n) $
under the Hermite semigroup has been studied by Byun \cite{By} . His
result is stated as follows.

Let $ \CH_t^\lambda $ be the Hermite-Bergman space defined to be the space
of all entire functions on $ \C^n $ for which
$$ \int_{\R^{2n}} |F(x+iy)|^2 U_t^\lambda(x,y) dx dy <\infty $$ 
where the weight function $ U_t $ is given by
$$ U_t(x,y) = c_n (\sinh(4\lambda t))^{-\frac{n}{2}}
e^{\lambda \tanh(2\lambda t)x^2}e^{-\lambda \coth(2\lambda t)y^2}.$$ 

\begin{thm}\label{thm=by}
The image
of $ L^2(\R^n) $ under the Hermite semigroup is precisely the space 
$ \CH_t^\lambda $ and $ e^{-tH(\lambda)} $ is a constant multiple of an
isometry between these two spaces. 
\end{thm}

The relation between the heat kernel
transform on $ \Gamma\bs\H $ and the Hermite semigroup is given in the
following proposition.

\begin{prop}\label{rel=hh}
Let $ f \in L^2(\R^n) $ and $ F = V_{k,\bj}(f) $ for $ {\bj} \in A_k.$ Then
$$ F*k_t(x,u,\xi) = c_\lambda e^{-t\lambda^2+i\lambda \xi}
e^{i\lambda ({\ba}\cdot x +\frac{1}{2}x\cdot u)}
\sum_{\bm \in \Z^n}e^{i\lambda x\cdot {\bm}}\tau_{-\ba}\left(e^{-tH(\lambda)}
\tau_{\ba}f\right)(u+{\bm})$$
where $ \lambda = 4\pi k , {\ba} = \frac{1}{2k}{\bj}, 
 \tau_{\ba}f(x)= f(x-{\ba}) $ and $  c_\lambda $ is a constant depending only
on $ \lambda $ and $ n.$
\end{prop}
\begin{proof} It follows from the definition of $ V_{k,\bj} $ and the
calculation (3.1) that $e^{-i\lambda \xi + t \lambda^2}F*k_t(x,u,\xi) $ is
given by

\begin{align*} 
&\int_{\R^{2n}}V_kf(x',u',0)e^{i\lambda {\ba}\cdot x}e^{i\frac{\lambda}{2}
(u\cdot x'-x\cdot u')}p_t^\lambda(x-x',u-u')dx'du'\\
&\quad =\sum_{\bm \in \Z^n}\int_{\R^{2n}}f(u'+{\bm})e^{i\lambda x'\cdot({\bm}+
\frac{1}{2}(u+u'+2{\ba}))}\cdot \\ 
&\quad \cdot e^{-i\frac{\lambda}{2}x\cdot u'}
p_t^\lambda(x-x',u-u')\ dx'\ du'\, .
\end{align*}

Using the explicit formula for $p_t^\lambda(x-x',u-u')$ the integral with
respect to $ dx' $ can be seen to be
\begin{align*}
& \int_{\R^{2n}}e^{i\lambda x'\cdot({\bm}+\frac{1}{2}(u+u'+2{\ba}))}
e^{-\frac{\lambda}{4}\coth(\lambda t)(x-x')^2}\  dx'\\
&\quad = e^{i\lambda x\cdot({\bm}+\frac{1}{2}(u+u'+2{\ba}))}\int_{\R^{2n}}
e^{-i\lambda x'\cdot({\bm}+\frac{1}{2}(u+u'+2{\ba}))}
e^{-\frac{\lambda}{4}\coth(\lambda t)x'^2}\  dx'\\
&\quad =   c_\lambda (\tanh(\lambda t))^{\frac{n}{2}}
e^{i\lambda x\cdot({\bm}+\frac{1}{2}(u+u'+2{\ba}))}
e^{-\lambda \tanh(\lambda t)({\bm}+\frac{1}{2}(u+u'+2{\ba}))^2}\ .
\end{align*}

Therefore,
\begin{align*}&  \int_{\R^{2n}}f(u'+{\bm})
e^{i\lambda x'\cdot({\bm}+\frac{1}{2}(u+u'+2{\ba}))}
e^{-i\frac{\lambda}{2}x\cdot u'} p_t^\lambda(x-x',u-u')\ dx'\ du'\\ 
&\quad  = c_\lambda (\sinh(2\lambda t))^{-\frac{n}{2}} \int_{\R^{2n}}f(u'+{\bm})
e^{i\lambda x\cdot({\bm}+\frac{1}{2}(u+u'+2{\ba}))}
e^{-i\frac{\lambda}{2}x\cdot u'} \\
& \quad \times e^{-\lambda \tanh(\lambda t)({\bm}+\frac{1}{2}(u+u'+2{\ba}))^2}
e^{-\frac{\lambda}{4}\coth(\lambda t)(u-u')^2}\ du'\ .\end{align*}

We change variables $ u' \rightarrow u'-{\ba}-{\bm} $, use the expression
for $ K_t^\lambda $ given in (\ref{Mehler})
and the integral above becomes 
\begin{align*}&  c_\lambda e^{i\lambda x\cdot {\bm}}e^{i\lambda x\cdot {\ba}}
e^{i\frac{\lambda}{2} x\cdot u }\int_{\R^n} f(u'-{\ba})
K_t^\lambda(u+{\ba}+{\bm},u')\ du'\\ 
&\quad = c_\lambda e^{i\lambda x\cdot {\bm}}e^{i\lambda x\cdot {\ba}}
e^{i\frac{\lambda}{2} x\cdot u }
\left(e^{-tH(\lambda)}\tau_{\ba}f\right)(u+{\ba}+{\bm})\, .
\end{align*}
This completes the proof of the proposition.
\end{proof}

We are ready for the main result in this article.

\begin{thm}\label{th=c} The map 
$$T_t: \CH_k \to \CB_{t,\Gamma}^\lambda, \ \ F\mapsto f*_\lambda p_t^\lambda$$
is an isometric isomorphism. 
\end{thm}

To prove the Theorem  we will establish 
the following  slightly more precise 
result.

\begin{thm}
An entire function $ F(z,w) $ belongs to $ \CB_{t,\bj}^{4\pi k} $ if and
only if $ F(z,w) = e^{t(4\pi k)^2}                                              
(V_{k,{\bj}}f)*k_t(z,w,0) $ for some $ f \in L^2(\R^n) .$
\end{thm}
\begin{proof} First note that the map is isometric by 
Proposition \ref{stand}
and Corollary \ref{corr}. It remains 
to verify surjectivity. 
\par For that let $ F \in 
\CB_{t,\bj}^{\lambda} $ with $ \lambda = 4\pi k.$ We have to show that there
exists $ f \in L^2(\R^n) $ such that 
$ V_{k,{\bj}}f*k_t(z,w,0)= e^{-t\lambda^2}F(z,w) $. 
To prove this we consider the function
$$ G(z,w) = e^{-i\lambda {\ba}\cdot z}e^{-i\frac{\lambda}{2} z\cdot w}F(z,w).$$
In view of the transformation properties of $ F $,  the function $ G $ becomes
$\frac{1}{2k}$-periodic in the $ x$-variables. Therefore, it admits an
expansion of the form
$$ G(z,w) = \sum_{\bm \in \Z^n} C_{\bm}(w) e^{i\lambda {\bm}\cdot z}$$
where $ C_{\bm} $ are the Fourier coefficients:
$$ C_{\bm}(w) = \int_{[0,\frac{1}{2k})^n}G(x,w)e^{-i\lambda {\bm}\cdot x}
\ dx.$$
The transformation properties of $ F $ lead to $ G(x,w-{\bm})= G(x,w)
e^{i\lambda {\bm}\cdot x}$ and hence  $ C_{\bm}(w-{\bm}) = C_{\bf{0}}(w).$
Thus, we obtain
$$ F(z,w) = e^{i\lambda {\ba}\cdot z}e^{i\frac{\lambda}{2} z\cdot w}
\sum_{\bm \in \Z^n} C_{\bf{0}}(w+{\bm})e^{i\lambda {\bm}\cdot z}.$$

We now show that $  C_{\bf{0}} $ belongs to the Hermite-Bergman space
$ \CH_t^\lambda$. For that 
we consider the finite integral: 

\begin{align*} \|F\|^2&= \int_{\C^{2n}/ \Z^{2n}}
|F(z,w)|^2 W_t^{-\lambda}(z,w)\ dx\ du\  dy\ dv \\ 
& = \int_{\C^{2n}/ \Z^{2n}} \left|e^{2i\lambda {\ba}\cdot z}e^{i\lambda z\cdot w}\right|\cdot 
\left|\sum_{\bm \in \Z^n} C_{\bf{0}}(w+{\bm})e^{i\lambda {\bm}\cdot z}\right|^2\cdot \\
&\quad \cdot  e^{-\lambda(u\cdot y - v\cdot x)} p_{2t}^{\lambda}(2y,2v)\ dx\ du\  dy\ dv \\
& = \sum_{\bm \in \Z^n} \int_{\C^{n}/ \Z^n}\int_{\R^n} e^{-2\lambda y \cdot(a+u+\bm)}
\left|C_{\bf{0}}(w+{\bm})\right|^2 \\ 
&\quad \cdot  p_{2t}^{\lambda}(2y,2v)\ du\  dy\ dv 
\end{align*}
We recall the explicit formula
for $p_t^\lambda$ and  use the fact 
$$ \int_{\R^n}e^{-2\lambda y\cdot u} e^{-\lambda \coth(2\lambda t)y^2} dy
= c_\lambda (\tanh(2\lambda t))^{\frac{n}{2}}e^{\lambda \tanh(2\lambda t)u^2}\, .$$
As a result 
$$\|F\|^2=\int_{\R^{2n}} |C_{\bf{0}}(w-{\ba})|^2 U_t^\lambda(u,v) \ du\ dv <\infty.$$
In view of Theorem \ref{thm=by}, 
there exists $ g \in L^2(\R^n) $ such that $C_{\bf{0}}(w) = 
e^{-tH(\lambda)}g(w+{\ba}).$ Let $ f = \tau_{-{\ba}}g $. Then 
Proposition \ref{rel=hh} implies 
$$ F(z,w) = e^{t\lambda^2} V_{k,\bj}f*k_t(z,w,0)$$ and 
this proves the theorem.
\end{proof}

\begin{rem} Our proof of the fact that the map 
$T_t\Big|_{\CH_k}$ is isometric (Proposition \ref{stand})
is rather robust and generalizes to all 
discrete subgroups $\Gamma<\H$. However this is not the case 
for our argument for the onto-ness.\end{rem}

\begin{rem}Recall the finite group 
${\bf F}_k$ and their two representations $\Pi_k$ and $\tilde \Pi_k$. 
Then $T_t\Big|_{\CH_k}$ intertwines 
$\Pi_k$ and $\tilde\Pi_k$. 
\end{rem}  

\begin{rem} The Weil-Brezin transforms $ V_{k,\bj} $ defined on
$ L^2(\R^n) $ has a natural extension to the Hermite-Bergman spaces
$ \CH_t^\lambda, \lambda = 4\pi k.$ Indeed, consider the operator $\tilde{
V}_{k,\bj}$ defined on $ \CH_t^\lambda $ as follows. For $ F \in 
\CH_t^\lambda $ we let
$$ \tilde{V}_{k,\bj}F(z,w,\zeta) = e^{i\lambda \zeta}e^{i\lambda {\ba}\cdot z}
e^{i\frac{\lambda}{2}z\cdot w}\sum_{\bm \in \Z^n} e^{i\lambda z\cdot {\bm}}
F(w+{\bm}).$$ Let us verify that the above series converges so that $
\tilde{V}_{k,\bj}F$ is well defined. As $ F \in  \CH_t^\lambda $ we have
$$ F(z) = \int_{\C^n} F(w) K_t^\lambda(z,\bar{w})U_t^\lambda(z) dz$$
since $ K_t^\lambda(z,\bar{w}) $ is the reproducing kernel. From the above we
get the estimate
$$ |F(z)| \leq C K_t^\lambda(z,\bar{z}) \leq C 
e^{-\lambda \tanh(2\lambda t)x^2}e^{\lambda \coth(2\lambda t)y^2}.$$ Therefore,
the series defining $\tilde{V}_{k,\bj}F(z,w,\zeta) $ converges uniformly
over compact subsets and defines an entire function. Moreover, we can check
that $ \tilde{V}_{k,\bj}F(z,w,0) $ belongs to $ \CB_{t,\bj}^\lambda.$
Thus $ \tilde{V}_{k,\bj} $ intertwines between the heat kernel transform
associated to the Hermite operator and the heat kernel transform on the
nilmanifold. More precisely, we have the following commutative diagram:

$$\begin{matrix}  L^2(\R^n) & \overset{V_{k,\bj}}\longrightarrow  & \CH_k \\
\downarrow & & \downarrow \\
\tau_{-\ba}e^{-tH(\lambda)}\tau_{\ba} &  & T_t\\
\downarrow & & \downarrow \\
\CH_t^\lambda & \overset{ \tilde{V}_{k,\bj}}\longrightarrow & \CB_{t,\bj}^{4\pi k}
\end{matrix}$$
\end{rem}

Let $\CB_{t,\Gamma}^0 $ be the Bergman space described in Theorem 4.1. Combining
Theorems 4.1 and 4.7 we obtain the following.

\begin{thm}
The image of $ L^2(\Gamma\bs\H) $ under $ T_t $ is the direct sum of all $ 
\CB_{t,\Gamma}^{4\pi k}, k \in \Z.$ More precisely,
$$ T_t(L^2(\Gamma\bs\H)) = \sum_{k=-\infty}^\infty e^{2t(4\pi k)^2}\CB_{t,\Gamma}^{4\pi k}.$$
\end{thm}

This is the analogue of Theorem 5.1 in \cite{KTX} for the heat kernel transform on the full
Heisenberg group. As in the case of the full Heisenberg group, the image can be written
as a direct sum of three weighted Bergman spaces if we allow the weight functions to take both
positive and negative values.

Let $ L^2_{+}(\Gamma\bs\H) =\oplus_{k=1}^\infty \CH_{-k} $ and  $ L^2_{-}(\Gamma\bs\H) =
\oplus_{k=1}^\infty \CH_{k}. $ Similarly define $ L^2_{+}(\T) $ and $ L^2_{-}(\T) $
where $ \T = \R / (\frac{1}{2}\Z)$ is the one dimesional torus. We let $ \CB_t^{+}(\C) $ (resp. $ \CB_t^{-}(\C) $)
stand for the image of $L^2_{+}(\T) $ (resp. $ L^2_{-}(\T) $) under the heat kernel transform
associated to the Laplacian on $ \T.$ These are weighted Bergman spaces that correspond to
the weight $ e^{-\frac{1}{2t}y^2} $ which are $ 1/2 $ periodic in the $ x-$ variable. We define
$\CB_t^{+}(\Gamma\bs\H_\C) $ and $ \CB_t^{-}(\Gamma\bs\H_\C) $ as follows. Let $ W_t^+ $ and $ W_t^-
$ be the weight functions that appeared in [7]. They are charactersied by the conditions
$$ \int_{\R} W_t^+(z,w,i\eta)e^{2\lambda \eta} \ d\eta = e^{2t\lambda^2} W_t^\lambda(z,w) $$
for all $ \lambda > 0 $ and
$$ \int_{\R} W_t^-(z,w,i\eta)e^{2\lambda \eta} \ d\eta = e^{2t\lambda^2} W_t^\lambda(z,w) $$
for all $ \lambda < 0.$  

We consider an exhaustion of $ \Gamma\bs\H $ defiend as follows. For each positive integer $ m$
let $ E_m = \{ z = x+iy \in \C^n : x \in [0,1)^n, |y| \leq m \} $ and $ E = [0,1) \times \R.$
We define $ K_m = E_m \times E_m \times E $ so that the union of all $ K_m $ as $ m $
varies over all positive integers is just $ \Gamma\bs\H.$ We let $ \CV^+(\Gamma\bs\H_\C) $ stand for
the space all entire functions $ F $ on $ \C^n \times \C^n \times \C $ such that $ F \in
L^2(K_m,|W_t^+|dg) $ for all $ m$;
$$ \lim_{m\rightarrow \infty} \int_{K_m}|F(g)|^2 W_t^+(g) \ dg < \infty;$$
and $ F(z,w,\cdot) \in \CB_t^{+}(\C) $ for all $ z,w \in \C^n.$ We equip  $ \CV^+(\Gamma\bs\H_\C) $
with the sesquilinear product
$$ (F,G)_+ = \lim_{m\rightarrow \infty}\int_{K_m} F(g)\bar{G}(g)W_t^+(g)\ dg.$$ We similarly
define $ \CV^-(\Gamma\bs\H_\C) .$
As in \cite{KTX}, Lemma 6.4, one shows that the above bracket is 
positive definite and we 
define $\CB_t^{+}(\Gamma\bs\H_\C) $ and $ \CB_t^{-}(\Gamma\bs\H_\C) 
$ as the completions of $ \CV^+(\Gamma\bs\H_\C) $ and $ \CV^-(\Gamma\bs\H_\C)$.
We also define
$ \CB_t^0(\Gamma\bs\H_\C) $ to be the space described in Theorem 4.1. With these definitions we
have the following result.

\begin{thm} 
$$ T_t(L^2(\Gamma\bs\H)) =  \CB_t^{-}(\Gamma\bs\H_\C) \oplus \CB_t^0(\Gamma\bs\H_\C) \oplus
\CB_t^{+}(\Gamma\bs\H_\C).$$
\end{thm} 

\begin{proof} The proof is similar to the case of the full Heisenberg group. So we will
be very sketchy. We refer to [7] for details (see Lemma 6.4 and Theorem 6.4). First we check
that the above sesquilinear product defines a pre-Hilbert space structure. The condition
$ F(z,w,\cdot) \in \CB_t^{+}(\C) $ shows that there exists a function $ G(z,w,\cdot) \in
L^2_+(\T) $ such that
$$  F(z,w,\zeta) = \int_{\R} G(z,w,s)q_t(\zeta-s) \ ds $$
where $ q_t$ is the Euclidean heat kernel. Therefore, upto a multiplicative constant
$$ \int_{[0,1/2)} F(z,w,\xi+i\eta) e^{-4\pi ik \xi} \ d\xi = e^{-4\pi k \eta} e^{-t(4\pi k)^2}
G^{-4\pi k}(z,w).$$ Since $ W_t^+ $ is independent of $ \xi $ we have
$$ \int_{K_m} |F(g)|^2 W_t^+(g) \ dg $$
$$ = \int_{E_m \times E_m}\int_{\R} 
\sum_{k=-\infty}^{-1}
|G^{-4\pi k}(z,w)|^2 e^{-2(4\pi k)\eta}e^{-2t(4\pi k)^2}W_t^+(z,w,i\eta)\ dz\ dw\ d\eta.$$
In view of the defining property of $ W_t^+ $ the above reduces to
$$ \int_{E_m \times E_m} \sum_{k=-\infty}^{-1}|G^{-4\pi k}(z,w)|^2 W_t^{-4\pi k}(z,w) \ dz\  dw.$$
This shows that the sesquilinear form is positive definite.

Consider the subspace $ L^2_{+,0}(\Gamma\bs\H) $ consisting of functions in
 $ L^2_{+}(\Gamma\bs\H) $ with only finitely many nonzero Fourier coefficients. Then it is
dense in $ L^2_{+}(\Gamma\bs\H)$. The proof will be complete if we can show that 
$$T_t( L^2_{+,0}(\Gamma\bs\H)) \subset  \CV^+(\Gamma\bs\H_\C) \subset 
T_t( L^2_{+}(\Gamma\bs\H)) $$ and the inclusions are isometric. 
Suppose $ F = f*k_t $ with $ f \in  L^2_{+,0}(\Gamma\bs\H). $ Then
$ f = \sum_{j=-N}^{-1} f_j, f_j \in \CH_j $ for some $ N $. Since $$ f_j*k_t(z,w,\zeta)
= e^{-t(4\pi j)^2}e^{4\pi ij\zeta}g_j*_{-4\pi j}p_t^{-4\pi j}(z,w)$$ where $ g_j(x,u)
= f_j(x,u,0)$
a simple  calculation shows that
$$ \lim_{m\rightarrow \infty} \int_{K_m}|F(g)|^2W_t^+(g)\ dg $$
$$=  \sum_{j=-N}^{-1} 
\int_{[0,1)^{2n}}\int_{\R^{2n}} |g_j*_{-4\pi j}p_t^{-4\pi j}(z,w)|^2 W_t^{-4\pi j}(z,w)dzdw$$
$$= C   \sum_{j=-N}^{-1} \int_{[0,1)^{2n}} |g_j(x,u)|^2 \ dx\  du = C \int_{\Gamma\bs\H}
|f(h)|^2 \ dh.$$ This proves that  $ T_t( L^2_{+,0}(\Gamma\bs\H))$ is isometrically embedded
into $ \CV^+(\Gamma\bs\H_\C)$.

To prove the inclusion $\CV^+(\Gamma\bs\H_\C) \subset T_t( L^2_{+}(\Gamma\bs\H)) $ let us
take $ F \in \CV^+(\Gamma\bs\H_\C) $ and consider the representation
$$  F(z,w,\zeta) = \int_{\R} G(z,w,s)q_t(\zeta-s) \ ds .$$ The argument used to prove
the positive definiteness of the inner product on  $\CV^+(\Gamma\bs\H_\C) $ shows that
$ G^{-4\pi k}(z,w) \in \CB_{t,\Gamma}^{-4\pi k} $ for $ -\infty < k \leq -1.$ We can write
$ G^{-4\pi k}(z,w) = g^{-4\pi k}*_{-4\pi k}p_t^{-4\pi k}(z,w) $ for some functions
$ g^{-4\pi k} .$ If we define
$$ f(x,u,\xi) = \sum_{-\infty}^{-1} e^{4\pi i k \xi} g^{-4\pi k}(x,u) $$ then it follows that
$ f*k_t(z,w,\zeta) = F(z,w,\zeta) $ as desired. This completes the proof.
\end{proof}

\end{document}